\documentclass[12pt,reqno]{amsart}

\usepackage{amsmath}
\usepackage{amssymb}
\usepackage{mathrsfs}
\usepackage{psfrag}
\usepackage{graphicx}

\begingroup
\newtheorem{theorem}{Theorem}[section]
\newtheorem{lemma}[theorem]{Lemma}
\newtheorem{proposition}[theorem]{Proposition}

\newtheorem{ass}{Assumption}
\endgroup

\theoremstyle{definition}
\begingroup
\newtheorem{definition}[theorem]{Definition}
\newtheorem{remark}[theorem]{Remark}

\endgroup
\numberwithin{equation}{section}

\newcommand{\bbibitem}{\bibitem}

\newcommand{\llabel}[1]{{\label{#1}}}

\newcommand{\NN}{{\mathbb N}}

\newcommand{\RR}{{\mathbb R}}

\newcommand{\h}{\nabla^2_x}
\newcommand{\n}{\nabla_x}
\newcommand{\e}{\varepsilon}
\newcommand{\D}{\delta}
\newcommand{\bx}{\bar{x}}

\newcommand{\ue}{u_\e}
\newcommand{\due}{\dot{u}_\e}

\title[]{Singular perturbations of finite dimensional gradient flows}

\author{Chiara Zanini}
\address[Chiara Zanini]{SISSA, Via Beirut 4, 34014 Trieste, Italy}
\email[Chiara Zanini]{zaninic@sissa.it}

\begin{document}

\footnotetext{Preprint SISSA 41/2006/M (July 2006).}

\noindent
\begin{abstract}
\noindent In this paper we give a description of the asymptotic behavior, as $\e\to 0$, of the $\e$-gradient flow in the finite dimensional case.

Under very general assumptions we prove that it converges to an evolution obtained by connecting some smooth branches of solutions to the equilibrium equation (slow dynamics) through some heteroclinic solutions of the gradient flow (fast dynamics).
\end{abstract}

\maketitle

\begin{section}{Introduction}\llabel{s:introd}

The study of quasistatic rate-independent evolutionary models may lead to consider gradient flow-type problems. 
Indeed, suppose that one wants to find a time-dependent function $t\mapsto u(t)$ satisfying
\begin{equation}\llabel{gf}
\n f(t,u(t))=0\,,
\end{equation}
where $f\colon[0,T]\times X\to \RR$ is a time-dependent energy functional, and $X$ is a given Banach space.
To this aim, it seems natural to study the perturbed problem
\begin{equation}\llabel{egf}
\e \dot{u}_\e(t) + \n f(t,u_\e(t)) = 0\,,
\end{equation}
which is indeed a gradient flow problem.
In order to obtain the existence of a quasistatic evolution, the limit as  $\e\to 0$ is analyzed.
The intention is to prove that under suitable assumptions on the functional $f$, the solutions $\ue$ converge to a limit function $u$ solving problem (\ref{gf}) and that this method selects the most interesting solutions~$u$ of (\ref{gf}) (see \cite{EM-03}, \cite{DM-Des-Mor-Mor-05-2}, \cite{Toa-Zan}).

In this paper we study a model case in which the main simplifying assumption is that the dimension of the space $X$ is finite.
We shall see that, under very general assumptions on $f$, the limit function $u(t)$ is a local minimum of $f(t,\cdot)$.
Moreover, it may admit some discontinuity times, while the approximating solutions $\ue(t)$ of the $\e$-gradient system (\ref{egf}) are at least continuous.

The first work on similar subjects was written by Efendiev and Mielke \cite{EM-03}, who add to the energy functional a dissipation term, which is crucial in the proof of the compactness of $\ue$. 
In our work, we do not have dissipative terms, but the assumptions on $f$ are stronger (see Section~\ref{s:set} below).

To be more specific, we consider a smooth energy function $f\colon [0,T]\times \RR^n\to \RR$ satisfying a suitable coerciveness condition (see Assumption~\ref{Hp1}). 
We suppose also that $\ue(0)\to u(0)$, $\n f(0,u(0))=0$, and that $\h f(0,u(0))$ is positive definite.
We will prove that $\ue$ converges, as $\e$ goes to zero, to a piecewise regular function $u\colon[0,T]\to \RR^n$, defined via the Implicit Function Theorem, such that $\n f(t,u(t))=0$ and $\h f(t,u(t))$ is positive definite on each continuity interval ${]t_{i-1},t_i[}$. 
It turns out that discontinuities $t_i$ of $u(t)$ are located at degenerate critical points of $f(t_i,\cdot)$, i.e., at points $x\in\RR^n$ where the Hessian matrix $\h f(t_i,x)$ possesses at least one zero eigenvalue.

To conclude this analysis we have to establish the connection between the limits $u(t_i^-)$ and $u(t_i^+)$.
This will be done by considering the fast dynamics, i.e., the dynamics governed by the system of differential equations
\begin{equation}\llabel{vintrod}
\dot{v}(s) = -\n f(t_i,v(s))\,.
\end{equation}

In a generic situation we may assume that $\h f(t_i,u(t_i))$ has exactly one zero eigenvalue, while the other eigenvalues are positive.
To discuss the behavior of (\ref{vintrod}), we are led to consider the systems of differential equations $\dot{v}(s)=-\n f(t,v(s))$ where $t$ is close to $t_i$ and plays the role of a parameter.
Under very general hypotheses, the following happens: before $t_i$ the vector field $\n f(t,\cdot)$ has two zeroes, a saddle and a node, at $t=t_i$ there is only one zero (the node and the saddle coalesce), and for $t>t_i$ these zeroes of the vector field no longer exist.
This corresponds to an abrupt change in the phase portrait as the parameter varies and is known in the literature as saddle-node bifurcation of codimension one (see \cite{GH-83}, \cite{Van}, \cite{HuWe-95}). 

In Section~\ref{s:set} we list the technical assumptions which permit to obtain the main result of the paper, Theorem~\ref{tm:main}.
Without entering all technical details, the setting obtained from our assumptions is the following one.
For every $t\in[0,T]$ there is a finite number of critical points $x\in\RR^n$ of $f(t,\cdot)$ and among them at most one is degenerate.
Moreover there exists only a finite number of pairs $(t,\xi)$ such that $\xi$ is a degenerate critical point of $f(t,\cdot)$.
On the degenerate critical points with only nonnegative eigenvalues, the Hessian matrix $\h f(t,\xi)$ has only one zero eigenvalue and satisfies two transversality conditions (see (b) and (c) in Assumption~\ref{Hp2}). 
Although we do not prove that Assumptions~\ref{Hp1}--\ref{Hp4} are generic in any technical sense, they cover a wide class of interesting examples.

If $\xi$ is a degenerate critical point of $f(t,\cdot)$ satisfying all conditions considered above, then we prove that there is a unique heteroclinic solution $v(s)$ of $\dot{v}(s)=-\n f(t,v(s))$ issuing from the degenerate critical point $\xi$, and we suppose that $v$ tends, as $s\to+\infty$, to a nondegenerate critical point $y$ of $f(t,\cdot)$, with $\h f(t,y)$ positive definite. 
The existence of such heteroclinic solution is standard. 
Since we have not been able to find the proof of uniqueness in the literature, we give the complete proof in Lemma~\ref{lm:evi}.

This analysis leads to a more precise construction of the function $u$ mentioned above.
Accordingly, the main result of this paper, Theorem~\ref{tm:main}, states that if $\ue(0)\to u(0)$, $\n f(0,u(0))=0$, and $\h f(0,u(0))$ is positive definite, and Assumptions~\ref{Hp1}--\ref{Hp4} are satisfied, then $\ue(t)$ converges to $u(t)$ uniformly on compact sets of $[0,T]\setminus\{t_1,\dots,t_{k-1}\}$, where $t_i$ are the discontinuity times for $u$. 
Moreover in a small neighborhood of $t_i$, a rescaled version of $\ue(t)$ converges to the heteroclinic solution $v(s)$, connecting $u(t_i^-)$ and~$u(t_i^+)$.
Finally, the graph of $\ue$ approaches the completion of the graph of $u$ obtained by using the heteroclinic trajectories.
\end{section}

\begin{section}{Setting of the problem} \llabel{s:set}

Throughout the paper, for fixed $T>0$, we make the following assumption:
\begin{ass}\llabel{Hp1}
$f\colon [0,T] \times \RR^n  \to \RR$ is a $C^3$-function satisfying the property
\begin{equation*}
\n f(t,x) \cdot x \geq c_0 \left|x\right|^2 - a_0\, , 
\end{equation*}
for some $a_0\geq 0$ and $c_0 > 0$,
\end{ass}
\noindent where $\n f=(f_{x_1},\dots,f_{x_n})$ denotes the gradient of $f$ with respect to its spatial variable $x\in\RR^n$.

We may deduce from this assumption that there exist two positive constants $M$ and $\tilde{c}$ (depending on $a_0$ and $c_0$), and a constant $\tilde{a}$ (depending also on $f$ and $T$) such that
\begin{equation}\llabel{finf}
f(t,x)\geq \tilde{c}|x|^2 - \tilde{a} \quad \mbox{ for every }|x|\geq M \mbox{ and every }t\in [0,T]\,.
\end{equation}
For given $t\in[0,T]$, we say that a point $x\in \RR^n$ is a critical point for $f(t,\cdot)$ if $\n f(t,x)=0$.
\begin{remark}\llabel{lm:ec}
Note that by Assumption~\ref{Hp1} all critical points for the function $f(t,\cdot)$ belong to the compact $\overline{B}$,
where $\overline{B} := \overline{B}(0, \sqrt{a_0c_0^{-1}})$ is the closed ball in $\RR^n$ centered at $0$ and with radius $\sqrt{a_0c_0^{-1}}$. Moreover, taking the minimum of $f(t,\cdot)$ in $\overline{B}$ it is immediate to get a critical point.
Hence, for every $t\in[0,T]$, critical points for $f(t,\cdot)$ exist and belong to $\overline{B}$. 
\end{remark}
We denote the set of zeroes to the gradient of $f$ by $\Gamma_f$, namely,
\begin{equation}\llabel{Gf}
\Gamma_f:=\{(t,x)\in [0,T]\times\RR^n:\n f(t,x)=0\}\,,
\end{equation}
and observe that $\Gamma_f\subset [0,T]\times\overline{B}$, by Remark~\ref{lm:ec}.

We recall that a critical point $\xi$ for $f(t,\cdot)$ is said to be degenerate if the kernel
of the Hessian matrix $\h f(t,\xi):=(f_{x_ix_j}(t,\xi))_{ij}$ is nontrivial, i.e., det$\h f(t,\xi)=0$.

In this paper, a particular interest will be brought on the set $Z_f$ of all pairs $(t,\xi)$ such that $\xi$ is a degenerate critical point for $f(t,\cdot)$, i.e.,
\begin{equation}\llabel{def:Z}
Z_f:= \{(t,\xi)\in\Gamma_f :  \det \h f(t,\xi)=0 \}.
\end{equation}
We make the following assumption.
\begin{ass} \llabel{Hp3} 
The number of all pairs $(t,\xi)$, such that $\xi$ is a degenerate critical point for $f(t,\cdot)$, is finite, i.e.,
\begin{equation*}
card (Z_f)=m<+\infty.
\end{equation*}
Moreover, let $\Pi\colon Z_f \to [0,T]$ denote the projection of $Z_f$ on the time-segment $[0,T]$, then we assume that $\Pi$ is injective and that $0,T\notin \Pi(Z_f)$.
\end{ass}
Throughout the paper we will focus on a particular class of degenerate critical points. 
More in detail, we make the following assumption.
\begin{ass}\llabel{Hp2}
For every $\tau\in[0,T]$ and for every degenerate critical point $\xi\in \RR^n$ for $f(\tau,\cdot)$, such that $\h f(\tau,\xi)$ is positive semidefinite, there exists $\ell\in \RR^n\setminus\{0\}$ such that the following conditions are satisfied:
\begin{itemize}
\item[(a)] $\ker \h f(\tau,\xi)={\rm span}(\ell) $;
\item[(b)] $ \n f_t(\tau,\xi)\cdot\ell\neq 0$, where $f_t(\tau,\xi)$ denotes the partial derivative of $f$ with respect to the time variable~$t$;
\item[(c)] $\sum_{i,j,k}f_{x_ix_jx_k}(\tau,\xi)\ell_i\ell_j\ell_k\neq 0$.
\end{itemize}
\end{ass}
Notice that condition (a) means that $0$ is a simple eigenvalue of $\h f(\tau,\xi)$ with eigenvector $\ell$, 
while the remaining $n-1$ eigenvalues are positive. 
Conditions (b) and (c) are known in the literature as transversality conditions (see, e.g., \cite{GH-83}). \begin{remark}\llabel{r:bec}
Let $(\tau,\xi)\in Z_f$, with $\h f(\tau,\xi)$ positive semidefinite.
An argument based on the Implicit Function Theorem (see, e.g., \cite{Van}),
implies that if $(\tau,\xi)$ satisfies Assumption~\ref{Hp2}, then there exists a smooth curve of solutions of $\n f(t(\lambda),x(\lambda))=0$, for $\lambda$ in a neighborhood of zero, with $(t(0),x(0))=(\tau,\xi)$.

More precisely, if conditions (b) and (c) have the same sign, then for every $t<\tau$ and near $\tau$ there are two solutions for the problem $\n f(t,x)=0$, while for $t$ near $\tau$ but $t>\tau$ there are no solutions. If conditions (b) and (c) have opposite sign, then the reverse is true.

Moreover,  the curve of zeroes passing through $(\tau,\xi)$ possesses a vertical tangent at $(\tau,\xi)$.
\end{remark}

\begin{remark}\llabel{r:Gf}
{}From our assumptions it turns out that $\Gamma_f$ is the union of a finite number of $C^2$-curves with end-points contained in $(\{0\}\times\RR^n)\cup(\{T\}\times \RR^n)$, see Figure~1.\ below for an example of the set $\Gamma_f$.
\end{remark}
\begin{figure}[h]
\begin{center}
\psfrag{t1}{$t_1$}
\psfrag{t2}{$t_2$}
\psfrag{t3}{$t_3$}
\psfrag{t4}{$t_4$}
\psfrag{m}{\scriptsize{ m}}
\psfrag{M}{ \scriptsize{M}}
\psfrag{0}{$0$}
\psfrag{T}{$T$}
\includegraphics[width=320pt]{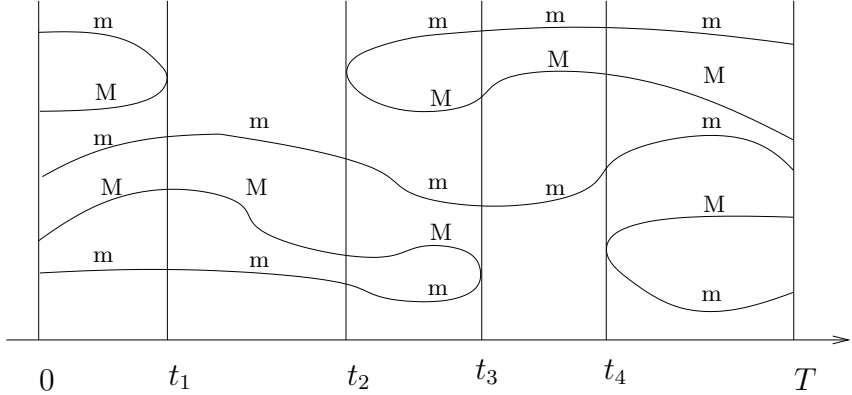}
\llabel{f:Gf}
\caption{\em An example for the set $\Gamma_f$: $m$ and $M$ stand for local minimum and maximum, respectively.}
\end{center}
\end{figure}
\begin{remark}\llabel{r:fincr}
The assumptions we made imply that for every $t\in [0,T]$ there exists a finite number of critical points for $f(t,\cdot)$. 
Indeed, if $t\notin Z_f$, then Assumption~\ref{Hp3} ensures that there are only nondegenerate critical points for $\n f(t,\cdot)$. Assumption~\ref{Hp1} implies that all critical points belong to the compact set $\overline{B}$, while by Assumption~\ref{Hp2} it follows that they are isolated.
On the other hand, by Assumption~\ref{Hp3}, at $t=\tau\in Z_f$ there is only one degenerate critical point $\xi$ for $f(\tau,\cdot)$.
\end{remark}

Let us freeze now a point $\tau\in \Pi(Z_f)$ and consider the autonomous system of differential equations in $\RR^n$ (depending on the single parameter $\tau$)
\begin{equation}\llabel{hv}
\dot{w}(s)=-\n f(\tau,w(s)).
\end{equation}
This is obviously a gradient system and, thanks to Assumption~\ref{Hp1}, (since positive semiorbits are bounded) we may apply the well known result that the $\omega$-limit set is contained into the set of equilibria of equation~(\ref{hv}) (see, e.g., \cite[Theorem 14.17]{HK}).
Moreover, since the equilibrium points are isolated (see Remark~\ref{r:fincr}), such an $\omega$-limit set is a single equilibrium point.

The following lemma ensures the existence of a unique heteroclinic solution $w$ issuing from $(\tau,\xi)\in Z_f$, while previous argument guarantees that $w$ has limit as $s\to+\infty$, and this limit is a (nondegenerate) critical point for $f(\tau,\cdot)$.
\begin{lemma}\llabel{lm:evi}
Suppose that Assumption~\ref{Hp1} and conditions (a) and (c) of Assumption~\ref{Hp2} are satisfied. 
Let $(\tau,\xi)$ be a point of $Z_f$ such that $\h f(\tau,\xi)$ is positive semidefinite.
Then there exists a unique (up to time-translations) solution of the problem
\begin{equation}\llabel{evi}
\begin{cases}
\dot{w}(s)=-\n f(\tau,w(s))\\
{\displaystyle \lim_{s\to-\infty} w(s)=\xi\,.}
\end{cases}
\end{equation}
\end{lemma}
\begin{proof} 
The proof is obtained by adapting a proof of the existence of the global center manifold, based on the Contraction Mapping Principle (see, e.g., \cite{Van}). 
The main difficulty is that usually, when the linearized part of a system of ordinary differential equations has some zero eigenvalue, there is, in general, existence of a heteroclinic solution, but not uniqueness (this is related to non-uniqueness of the local center manifold, see, e.g., \cite{GH-83}, \cite[\S 1.4]{Van}). Here the uniqueness is obtained thanks to the particular conditions (a) and (c) of Assumption~\ref{Hp2}.

During the proof, we will use the following notation: $g(x):=f(\tau,x)$, for every $x\in\RR^n$. 
To simplify further the formulation we make a number of preliminary transformations: a translation to take $\xi$ to the origin, and a linear transformation to bring $\nabla^2 g(0)$ in a diagonal form where the first eigenvalue is zero with eigenvector $e_1=(1,0,\dots,0)$.
Therefore we are reduced to the following hypotheses: 
\begin{equation}\llabel{hpg}
\nabla g(0)=0, \qquad
\nabla^2 g(0)=\left(\begin{array}{cc}
0&0\\
0&A
\end{array}\right)\,,\qquad \mbox{and}\qquad g_{x_1x_1x_1}(0)\neq 0\,,
\end{equation}
where $A$ is an $(n-1)\times(n-1)$ diagonal and invertible matrix. 
Moreover, by our assumption, the diagonal entries of $A$ are all positive real numbers.
In order to simplify the notation, we also suppose that $\frac{1}{2} g_{x_1x_1x_1}(0)=1$.

In order to obtain existence and uniqueness of the problem
\begin{equation}\llabel{evi1}
\begin{cases}
\dot{w}(t)=-\nabla g(w(t)),\\
{\displaystyle \lim_{t\to-\infty}w(t)=0}\,,
\end{cases}
\end{equation}
we apply the Contraction Mapping Theorem.
More precisely, for every $x\in\RR^n$ we use the following decomposition $x=(x_1,\bx)$, with $\bx:=(x_2,\dots,x_n)\in\RR^{n-1}$ and consider the space $Y$ of all functions $y\colon (-\infty,0]\to\RR^n$, $y(t)=(y_1(t),\bar{y}(t))$, such that 
\begin{equation*}
\|y_1\|_{Y_1}:=\sup_{ t\leq 0}|(t-1)y_1(t)|<\infty\,,\quad\mbox{and}\quad 
\|\bar{y}\|_{\overline{Y}}:=\sup_{t\leq 0}|(t-1)^2\bar{y}(t)|<\infty\,,
\end{equation*}
endowed with the norm
\begin{equation*}
\|y\|_Y:= \|y_1\|_{Y_1}+\|\bar{y}\|_{\overline{Y}}\,.
\end{equation*}
For every $x\in\RR^n$ let $\nabla g(x)=(D_1g(x),\bar{D}g(x))$.
Using now the Taylor expansion for $x$ in a neighborhood of $0\in\RR^n$, we get
\begin{equation}\llabel{Dg}
\begin{array}{c}
{\displaystyle D_1g(x)= 
x_1^2+x_1b\cdot\bx+\varphi(\bx,\bx)+o(|x|^2)}\\
{\displaystyle \bar{D}g(x)=A\bx + x_1^2b+x_1B\bx+\Phi(\bx,\bx)+o(|x|^2)}\,,
\end{array}
\end{equation}
where $b$ is a suitable vector in $\RR^{n-1}$, $\varphi\colon\RR^{n-1}\times\RR^{n-1}\to \RR$ and $\Phi\colon\RR^{n-1}\times\RR^{n-1}\to\RR^{n-1}$ are bilinear symmetric forms (whose coefficients depend on the third derivative of $g$ at the origin), $A$ is the matrix which appears in (\ref{hpg}), and $B$ is a $(n-1)\times(n-1)$ matrix whose entries depend on the third derivative of $g$ at $0$. 
Moreover, for every $y\in Y$ we define 
\begin{equation}\llabel{hhh}
\begin{array}{c}
{\displaystyle h_1(t):=-D_1g(y(t))+y_1(t)^2}\\
{\displaystyle \bar{h}(t):= -\bar{D}g(y(t))+A\bar{y}(t)}\,,
\end{array} 
\end{equation}
and observe that due to (\ref{Dg}) the asymptotic behavior at $-\infty$ is $h_1(t)\sim (t-1)^{-2}$, and $\bar{h}(t)\sim (t-1)^{-2}$, respectively.
For every $h=(h_1,\bar{h})$ satisfying these estimates, let us consider the function $x=(x_1,\bar{x})$ obtained by solving the following two problems depending on a parameter $\e>0$, which will be fixed later.
\begin{equation}\llabel{I'}
\begin{cases}
\dot{x}_1+\e x_1^2=\e h_1(t) &\mbox{ on }(-\infty,0],\\
{\displaystyle x_1(0)=-\frac{1}{\e}\,},
\end{cases}
\end{equation}
and 
\begin{equation}\llabel{II'}
\begin{cases}
\dot{\bx}+A\bx=\e\bar{h}(t) &\mbox{ on }(-\infty,0],\\
{\displaystyle\lim_{t\to-\infty}\bx(t)=0\,.}
\end{cases}
\end{equation}
We shall prove that problem~(\ref{II'}) has a unique solution with $\|\bx\|_{\overline{Y}}$ finite, and that for $\e$ sufficiently small the solution of problem~(\ref{I'}) does exist and satisfies
\begin{equation}\llabel{limx1}
\lim_{t\to-\infty}x_1(t)=0\,.
\end{equation}
Note that, if $h=(h_1,\bar{h})$ is defined by (\ref{hhh}), and if $x=y$, then $w:=\e x$ solves problem~(\ref{evi1}).

{}From the variation of constant formula it follows that the unique solution of problem (\ref{II'}) is
\begin{equation}\llabel{bxx}
\bar{x}(t) = \e \int_{-\infty}^t e^{-A(t-s)}\bar{h}(s)\, ds\,.
\end{equation}
Moreover $\|\bar{x}\|_{\overline{Y}}<\infty$.

Let us discuss now the existence of a solution of problem~(\ref{I'}).
We claim that for $\e$ sufficiently small there exists a solution defined on $(-\infty,0]$.
This is done using differential inequalities, since there exists a positive constant $M$ such that $|h_1(t)|\leq M/(t-1)^2$. 
More precisely, we are reduced to study two auxiliary problems, with the same condition in zero:
\begin{equation}\llabel{Iaux}
\begin{cases}
{\displaystyle \dot{x}_1 + \e x_1^2 =  \frac{M}{(t-1)^2}} &\mbox{ on }(-\infty,0]\\
{\displaystyle x_1(0)=-\frac{1}{\e}}\,.
\end{cases}
\end{equation} 
and 
\begin{equation}\llabel{IIaux}
\begin{cases}
{\displaystyle \dot{x}_1 + \e x_1^2 = - \frac{M}{(t-1)^2}} &\mbox{ on }(-\infty,0]\\
{\displaystyle x_1(0)=-\frac{1}{\e}}\,.
\end{cases}
\end{equation} 
For our purposes, it is sufficient to prove that both solutions of the auxiliary problems (\ref{Iaux}) and (\ref{IIaux}), tend to zero as $t\to-\infty$, for $\e$ sufficienlty small (depending on $M$).

For the existence, we observe that the equation considered in (\ref{Iaux}) (and in (\ref{IIaux}), respectively) is a particular case of the Riccati equation (see, e.g., \cite{MR}). 
Putting $x_1(t)= u(t)/(\e(t-1))$ (in both cases) we obtain an equation in the unknown $u(t)$ which is of the first order and can be solved by separation of variables. 
Moreover, if $\e$ is sufficiently small, then the solution $u(t)$ related to problem (\ref{Iaux}) is defined on $(-\infty,0]$ and is bounded at $-\infty$. 
On the other hand, the solution $u(t)$ related to problem (\ref{IIaux}) is defined on $(-\infty,0]$ and is bounded at $-\infty$, for every $\e>0$. 
Hence, for $\e$ small enough, we obtain an upper (lower) function solving the auxiliary problem (\ref{Iaux}) (or (\ref{IIaux}), respectively), and tending to zero as $t\to-\infty$.

Using differential inequalities (see, e.g., \cite[Theorem 6.1]{Hale}), we deduce that there exists $\e_0=\e_0(M)$ such that problem (\ref{I'}) admits a unique solution satisfying also the limit condition (\ref{limx1}), for every $\e<\e_0$.
Moreover, it is immediate to prove that the asymptotic behavior of $x_1(t)$ at $-\infty$ is like $(t-1)^{-1}$.

Finally, for $\e<\e_0$ we can define the map $\Gamma\colon Y\to Y$ by setting
\begin{equation}\llabel{Gamma}
\Gamma(y)(t) := (x_1(t),\bx(t))\,, 
\end{equation}
where $x_1(t)$ is the solution of (\ref{I'})--(\ref{limx1}), and $\bx(t)$ is given by (\ref{bxx}), with $h_1$ and $\bar{h}$ defined by (\ref{hhh}). 
Obviously, $\Gamma(y)$ belongs to $Y$, while it remains to prove that the map $y\mapsto \Gamma(y)$ is a strict contraction, for $\e$ sufficiently small.

Let us begin with the first component of $\Gamma(y)$.
For every $y\in Y$, let $h(t)=(h_1(t),\bar{h}(t))$ be defined as in (\ref{hhh}), and let us pass from $t$ to $-t$.
Then, there exists $H_1\in L^\infty(0,+\infty)$ such that $-h_1(-t)=\frac{H_1(t)}{(1+t)^2}\,$ for every $t\geq 0$. Let $x_1(t)$ be the solution to the following problem 
\begin{equation*}
\begin{cases}
{\displaystyle \dot{x}_1-\e x_1=\e\frac{H_1(t)}{(1+t)^2}} &\mbox{ on }[0,+\infty)\\
{\displaystyle x_1(0)=-\frac{1}{\e}\,.}
\end{cases}
\end{equation*}
In the same way, starting from $y^*\in Y$, we define $h^*=(h^*_1,\bar{h}^*)$, $H^*_1\in L^\infty(0,+\infty)$, and $x_1^*$ as the solution of an analogous problem having $H^*_1$ in the right-hand side, instead of $H_1$.
Put $x_1(t)=\frac{u(t)}{\e(1+t)}$, so that $u(t)$ solves the problem
\begin{equation*}
\begin{cases}
(1+t)\dot{u}=\e^2H_1(t)+u^2+u&\mbox{ on }[0,+\infty)\\
u(0)=-1\,.
\end{cases}
\end{equation*}
By this choice of the initial datum, we deduce that 
\begin{equation}\llabel{Mu}
|u(t)+1|\leq \e^2M\,,
\end{equation}
for every $t\geq 0$, being $M$ an upper bound for the $L^\infty$-norm of $H_1$.
Arguing in the same manner for $x_1^*$, we define $u^*$.
We want to prove that
\begin{equation}\llabel{1uuu}
|u(t)-u^*(t)|\leq \e^2C\|H_1-H^*_1\|_\infty\,,
\end{equation}
for every $t\geq 0$, so that, passing from $t$ to $-t$ and not renaming the solutions $x_1$ and $x_1^*$, we will get
\begin{equation*}
|x_1(t)-x_1^*(t)|\leq \frac{\e}{|t-1|} C\|y-y^*\|_Y\,\quad\mbox{ for every }t\leq 0\,,
\end{equation*}
where we used the inequality $\|H_1-H^*_1\|_\infty \leq C\|y-y^*\|_Y$, which follows from (\ref{hhh}).
We will obtain that
\begin{equation}\llabel{1xxx}
\|x_1-x_1^*\|_{Y_1}\leq \frac{1}{2}\|y-y^*\|_Y\,,
\end{equation}
having supposed that $\e C<\frac{1}{2}\,$.

Therefore, we are reduced to prove (\ref{1uuu}). 
Let $z(t):=u(t)-u^*(t)$, and $\alpha(t):=-u(t)-u^*(t)>0$.
Then $z(t)$ solves the problem
\begin{equation*}
\begin{cases}
(1+t)\dot{z} = \e^2(H_1(t)-H_1^*(t)) -\alpha(t)z + z &\mbox{ on }[0,+\infty)\\
z(0)=0\,.
\end{cases}
\end{equation*}
By the variation of constant method, we deduce that the solution $z(t)$ can be represented by the following formula:
\begin{equation*}
z(t)=\e^2\int_0^t \frac{H_1(s)-H_1^*(s)}{1+s}\, e^{-\int_s^t\frac{\alpha(\sigma)-1}{1+\sigma}\,d\sigma}\,ds\,.
\end{equation*}
Due to (\ref{Mu}), we obtain that $u(t)< -\frac{3}{4}$ for $\e$ sufficiently small, 
and the same is true for $u^*(t)$. 
Hence, $\alpha(t)> \frac{3}{2}$, for $\e$ sufficiently small.
Then
\begin{equation*}
\begin{array}{c}
{\displaystyle |z(t)|\leq \e^2\|H_1-H_1^*\|_\infty \Big|\int_0^t\frac{1}{1+s}\, e^{-\frac{1}{2}\int_s^t\frac{d\sigma}{1+\sigma}}\,ds\,\Big|=}\\
{\displaystyle = \e^2\|H_1-H_1^*\|_\infty \frac{1}{(1+t)^{\frac{1}{2}}}\Big|\int_0^t (1+s)^{-\frac{1}{2}}\, ds\Big|\leq 2\e^2\|H_1-H_1^*\|_\infty }\,.
\end{array}
\end{equation*}
This last estimate gives (\ref{1uuu}), and therefore (\ref{1xxx}) is proved.

Let us consider now the second component of $\Gamma(y)$.
Let $\Gamma(y)(t)=(x_1(t),\bx(t))$ and $\Gamma(y^*)(t)=(x_1^*(t),\bx^*(t))$.
Therefore,
\begin{equation*}
\bx(t)-\bx^*(t) = \e\int_{-\infty}^t (\bar{h}(s)-\bar{h}^*(s))e^{-A(t-s)}\,ds\,.
\end{equation*}
Hence, using the fact that, by (\ref{hhh}), $|\bar{h}(s)-\bar{h}^*(s)|\leq C \|y-y^*\|_Y (s-1)^{-2}$, we get
\begin{equation*}
(t-1)^2|\bx(t)-\bx^*(t)| \leq C\e(t-1)^2\int_{-\infty}^t\frac{e^{-A(t-s)}}{(s-1)^{2}}\,ds\|y-y^*\|_Y \,.
\end{equation*}

Since
\begin{equation*}
\sup_{t\leq 0} \Big((t-1)^2\int_{-\infty}^t\frac{e^{-A(t-s)}}{(s-1)^{2}}\,ds\Big) < +\infty\,,
\end{equation*}
we deduce that there exists a positive constant $C^*$ such that $(t-1)^2|\bar{x}(t)-\bar{x}^*(t)|\leq \e C^*\|y-y^*\|_Y$, i.e., 
\begin{equation}\llabel{bxxx}
\|\bx-\bx^*\|_{\bar{Y}}\leq \frac{1}{2}\|y-y^*\|_Y\,,
\end{equation}
having supposed that $\e C^*<\frac{1}{2}$.

This estimate, together with (\ref{1xxx}), guarantees that the inequality
\begin{equation*}
\|\Gamma(y)-\Gamma(y^*)\|_Y\leq \frac{1}{2}  \|y-y^*\|_Y
\end{equation*}
holds true, and this concludes the proof.
\end{proof}
Throughout the paper, we make the following assumption.
\begin{ass}\llabel{Hp4}
For every $(\tau,\xi)\in [0,T]\times\RR^n$ such that $\xi$ is a degenerate critical point for $f(\tau,\cdot)$, satisfying the assumptions of Lemma~\ref{lm:evi}, let $w$ be the unique solution of (\ref{evi}) corresponding to $\tau$ and $\xi$. 
Let $w_\infty:=\lim_{s\to +\infty}w(s)$, then we assume that
\begin{equation}\llabel{ip4}
\h f(\tau,w_\infty)\quad\mbox{ is positive definite\,}.
\end{equation}
\end{ass}

\end{section}

\begin{section}{Preliminary results}\llabel{s:pres}

Starting from a suitable point $(\bar{t},\bx)\in{[0,T[}\times\RR^n$ we prove in the next lemma the existence of a maximal interval ${[\bar{t},\hat{t}[}$, and of a regular function $u$, defined on ${[\bar{t},\hat{t}[}$, such that $u(t)$ is a critical point for $f(t,\cdot)$ for every $t\in {[\bar{t},\hat{t}[}$.
\begin{lemma}\llabel{lm:u-ext}
Let $0\leq \bar{t}<T$, and let $\bx\in \RR^n$ be such that $\n f(\bar{t},\bx) = 0$ and $\h f(\bar{t},\bx)$ is positive definite. 
Suppose that Assumptions~\ref{Hp1} and \ref{Hp3} are satisfied.
Then there exist a maximal interval of existence ${[\bar{t},\hat{t}[}$, and a function $u\colon{[\bar{t},\hat{t}[}\to \RR^n$ of class $C^2$, such that $u(\bar{t})=\bx$ and $\n f(t,u(t))=0$ for every $t\in {[\bar{t},\hat{t}[}$.
Moreover, either $\hat{t}=T$ or $\hat{t}$ belongs to $\Pi(Z_f)$ (defined in Assumption~\ref{Hp3}).
\end{lemma}
\begin{proof}
The Implicit Function Theorem ensures that there are a maximal interval of existence ${[\bar{t},\hat{t}[}$ and a function 
$u\colon{[\bar{t},\hat{t}[}\to \RR^n$ of class $C^2$ such that 
\begin{equation*}
\n f(t,u(t))=0\qquad\mbox{and}\qquad \det\h f(t,u(t))> 0\quad\mbox{ on }{[\bar{t},\hat{t}[}.
\end{equation*}
The next step is to prove that $u(t)$ has limit as $t$ approaches $\hat{t}$ to the left. 
This is trivial if $\hat{t}=T$.
For the case $\hat{t}<T$, we introduce the following auxiliary result that will be proved later.
\begin{lemma}\llabel{lm:omlim}
Under the same assumptions of Lemma~\ref{lm:u-ext}, let us define the following set
\begin{equation}\llabel{K}
K:= \{x\in \RR^n| \exists s_k \nearrow \hat{t}: u(s_k)\to x\}.
\end{equation}
Then $K$ is a compact and connected set, composed only of critical points of $f(\hat{t},\cdot)$. Moreover, if $\hat{t}<T$ then $\det \h f(\hat{t},x)=0$ for any $x\in K$.
\end{lemma}
\noindent{\em Proof of Lemma~\ref{lm:u-ext} (continued).}
Let us suppose that Lemma~\ref{lm:omlim} is true, and let us prove that
\begin{equation}\llabel{limu1}
\lim_{t\to \hat{t}^-} u(t)
\end{equation}
does exist.
Indeed, let $K$ be the nonempty set defined by (\ref{K}). 
We need to show that $K$ reduces to just one point. 
Assume by contradiction that the limit (\ref{limu1}) does not exist. 
Then there are at least two sequences $s^i_k \nearrow \hat{t}$, $i=1,2$ and two distinct 
points $w_1,w_2\in \RR^n$ such that $u(s^i_k)\to w_i$, $i=1,2$. 
But $K$ is a connected set, thus there exists a continuous path $\gamma\colon [0,1]\to \RR^n$, connecting $w_1$ to $w_2$, 
such that $\gamma([0,1])\subset K$. 
The contradiction comes from the fact that by Lemma~\ref{lm:omlim} $x\in K$ implies 
$(\hat{t},x)\in Z_f$, which is finite by Assumption~\ref{Hp3}.

Finally, if $\hat{t}<T$, then by Lemma~\ref{lm:omlim} $\det \h f(\hat{t},x)=0$ for any $x\in K$, while by continuity $\n f(\hat{t},x)=0$, i.e., every $x\in K$ is a degenerate critical point for $f(\hat{t},\cdot)$. 
Hence by Assumption~\ref{Hp3}, $\hat{t}\in \Pi(Z_f)$, and this concludes the proof of Lemma~\ref{lm:u-ext}.
\end{proof}

\begin{proof}[Proof of Lemma \ref{lm:omlim}.]
We begin with compactness. By definition the set $K$ is closed, while Assumption~\ref{Hp1} guarantees that it is bounded (see Remark~\ref{lm:ec}).

We continue by proving that $K$ is connected. This can be done in two steps.
The first one consists into prove that for any neighborhood $U$ of the set $K$ there exists $k>0$ 
such that $u(s)\in U$ for any $s\in V_k:={[\hat{t} - \frac 1 k,\hat{t} [}$, that is, in other words, 
$u(s)$ converges to $K$ whenever $s\to \hat{t}$. 
This can be done arguing by contradiction and using again Assumption~\ref{Hp1}.
The second step consists in taking two closed and disjoint sets $A$ and $B$ and assuming by contradiction $B\cap K=K\setminus A$, and that distance$(A\cap K,B\cap K)$ is positive. 
Then the first step gives the contradiction.
These two arguments are standard and we omit the details of them.

Last, let $x\in K$ and assume $\hat{t}<T$. Then by definition there exists $s_k\nearrow\hat{t}$ such that $u(s_k)$ approaches $x$ as $k\to\infty$.
By continuity, $\det \h f(s_k,u(s_k))$ tends to $\det \h f(\hat{t},x)$, and, moreover, $\n f(\hat{t},x)=0$.
If $\det \h f(\hat{t},x)\neq 0$, then the Implicit Function Theorem could be applied, a contradiction with the definition of $\hat{t}$. This concludes the proof.
\end{proof}
Starting from $\bar{t}=0$ and from a suitable point $y_0\in\RR^n$, we may repeatedly apply Lemma~\ref{lm:u-ext} and Lemma~\ref{lm:evi} obtaining the result stated in the following proposition.
\begin{proposition}\llabel{p:u}
Suppose that Assumptions~\ref{Hp1}--\ref{Hp4} are satisfied.
Let $y_0$ be such that $\n f(0,y_0)=0$ and $\h f(0,y_0)$ is positive definite.
Then there exist a unique (and finite) family of times $0=t_0<t_1<\dots<t_{k-1}<t_k=T$ and a unique family of functions $u_i\colon{[t_{i-1},t_i[}\to \RR^n$ of class $C^2$, for $i=1,\dots,k$, and a unique (up to time-translations) family of functions $v_i\colon\RR\to\RR^n$ of class $C^2$, $i=1,\dots,k-1$, such that
\begin{itemize}
\item[(1)] $u_1(0)=y_0$,
\item[(2)] for every $t\in {[t_{i-1},t_i[}$, $\n f(t,u_i(t))=0$ and $\h f(t,u_i(t))$ is positive definite,
\item[(3)] for every $i=1,\dots, k$, there exists $x_i:=\lim_{s\to t_i^-}u_i(t)$, and for every $i=1,\dots,k-1$, $(t_i,x_i)\in Z_f$, $\h f(t_i,x_i)$ is positive semidefinite and conditions (b) and (c) of Assumption~\ref{Hp2} have the same sign,
\item[(4)] for every $i=1,\dots,k-1$, function $v_i(s)$ solves 
\begin{equation}\llabel{vi}
\dot{v}_i(s)=-\n f(t_i,v_i(s))
\end{equation}
and satisfies 
\begin{equation}\llabel{vii}
\lim_{s\to-\infty}v_i(s) = \lim_{t\to t_i^-} u_i(t)\qquad \qquad\lim_{s\to +\infty} v_i(s)=u_{i+1}(t_i).
\end{equation}
\end{itemize}
\end{proposition}
\begin{proof}
We apply Lemma~\ref{lm:u-ext} with $(\bar{t},\bx)=(0,y_0)$ obtaining the existence of $\hat{t}=:t_1$, and of a function $u_1\colon{[0,t_1[}\to\RR^n$ of class $C^2$ such that $u_1(0)=y_0$, $\n f(t,u_1(t))=0$ on ${[0,t_1[}$, and $\h f(t,u_1(t))$ is positive definite on ${[0,t_1[}$. 
This proves conditions (1) and (2) restricted to ${[0,t_1[}$.

Arguing as in the proof of Lemma~\ref{lm:u-ext} (using Lemma~\ref{lm:omlim}), we deduce that there exists $x_1:=\lim_{t\to t_1^-}u_1(t)$ and all eigenvalues of $\h f(t_1,x_1)$ are nonnegative.
Moreover, since for every $t<t_1$ the function $u_1(t)$ solves the problem $\n f(t,x)=0$, then it follows from Remark~\ref{r:bec} that the transversality conditions (b) and (c) of Assumption~\ref{Hp2} have the same sign.
Indeed, if on the contrary they had the opposite sign, then there should be no solutions for the problem $\n f(t,x)=0$ for $t$ belonging to a left neighborhood of $t_1$. 
Thus condition (3) (restricted to ${[0,t_1[}$) is satisfied.

By Lemma~\ref{lm:evi} there exists a unique (up to time-translations) heteroclinic solution $v_1$ issuing from $x_1$.
In addition, as $s\to+\infty$, $v_1(s)$ tends to a critical point $y_1$ for $f(t_1,\cdot)$, and by Assumption~\ref{Hp4}, $\h f(t_1,y_1)$ is positive definite, so that condition (4) for $i=1$ is satisfied.

Next, if $t_1<T$, we apply Lemma~\ref{lm:u-ext} with $(\bar{t},\bx)=(t_1,y_1)$ and repeat the previous arguments.
\end{proof}
\begin{definition}\llabel{d:u}
Suppose that Assumptions~\ref{Hp1}--\ref{Hp4} are satisfied.
For fixed $y_0\in \RR^n$ such that $\n f(0,y_0)=0$ and $\h f(0,y_0)$ is positive definite, 
let $u_i$, $i=1,\dots,k$ be the functions obtained in Proposition~\ref{p:u}.
We thus define the $C^2$-piecewise function $u\colon [0,T]\to \RR^n$, such that $u(0)=y_0$, by
\begin{equation*}
u_{\big|{[t_{i-1},t_{i})}}:=u_i,\qquad \mbox{ for every }i=1,\dots,k\,.
\end{equation*}
Hence, $u$ is discontinuous at $t_1<\dots<t_{k-1}$ and satisfies 
\begin{equation}\llabel{du}
\n f(t,u(t))=0 \quad \mbox{for every }t\in{[t_{i-1},t_{i})},\quad u(t_{i-1})=y_{i-1}\quad\mbox{and}\quad\lim_{t\to t_{i}^-}u(t)=x_i\,,
\end{equation}
for every $i=1,\dots, k$ (cfr. Figure~2.).
\end{definition}
\begin{figure}[ht]\llabel{fig:du}
\begin{center}
 \psfrag{y0}{$y_0$}
 \psfrag{y1}{$y_1$}
 \psfrag{y2}{$y_2$}
 \psfrag{y3}{$y_{3}$}
 \psfrag{x1}{$x_1$}
 \psfrag{x2}{$x_2$}
 \psfrag{x3}{$x_{3}$}
 \psfrag{x4}{$x_4$}
 \psfrag{t1}{$t_1$}
 \psfrag{t2}{$t_2$}
 \psfrag{t3}{$t_{3}$}
 \psfrag{u1}{$u_1$}
 \psfrag{u2}{$u_2$}
 \psfrag{u3}{$u_{3}$}
 \psfrag{u4}{$u_4$}
 \psfrag{0}{$0=t_0$}
 \psfrag{T}{$t_4=T$}
\psfrag{dots}{$\dots$}
\includegraphics[width=300pt]{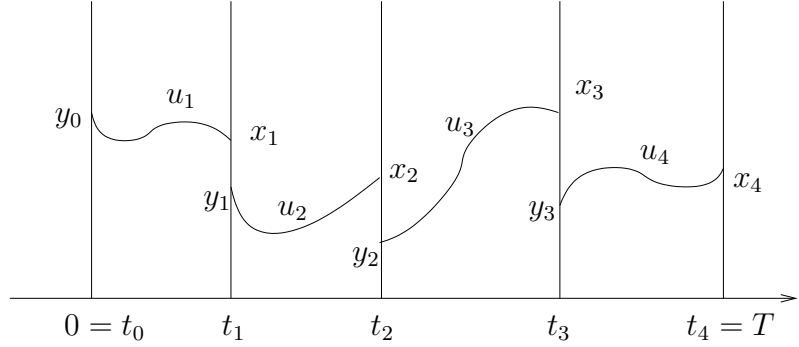}
\caption{\em The $C^2$-piecewise function $u$, expressed in terms of the functions $u_i$ defined in Proposition~\ref{p:u}, when $k=4$.}
\end{center}
\end{figure}
At $(t_i,x_i)$ by Assumption~\ref{Hp2} the Hessian matrix $\h f(t_i,x_i)$ has one zero eigenvalue while by construction the remaining $n-1$ eigenvalues are positive, for $i=1,\dots,k-1$. 
By Remark~\ref{r:bec} there exist $r_i>0$ and $R_i>0$ such that the following conditions hold true (see also Figure~3.).
\begin{itemize}
\item[$(C_l)$] There are two regular branches of solutions of $\n f(t,x)=0$ for $t\in[t_i-r_i,t_i]$, $i=1,\dots,k-1$. Moreover, if ker$\n f(t_i,x_i)={\rm span}(\ell_i)$, then the two branches have common (vertical) tangent $(0,\ell_i)$ at $(t_i,x_i)$, $i=1,\dots,k-1$;
\item[$(C_r)$] we have
\begin{equation}\llabel{ehi}
|\n f(t,x)|>0, \qquad\mbox{ on } {]t_i,t_i+r_i]}\times \overline{B}(x_i,R_i)
\end{equation}
for every $i=1,\dots,k-1$. 
\end{itemize}
By Definition~\ref{d:u} and condition $(C_l)$ one of these two regular branches of solutions has graph contained in $\{(t,u(t))|t\in[0,T]\setminus\{t_1,\dots,t_k\}\}$. 
Throughout the paper, the other one 
will be called $\bar{u}(t)$. 
It follows that $u$ and $\bar{u}$ have the same limit as $t\to t_i^-$.
\begin{figure}[ht]\llabel{f:ubaru}
\begin{center}
\psfrag{df}{$|\nabla_x f(t,x)|>0$}
\psfrag{ut}{$u(t)$}
\psfrag{but}{$\bar{u}(t)$}
\psfrag{-tiri}{$t_i-r_i$}
\psfrag{ti}{$t_i$}
\psfrag{tiri}{$t_i+r_i$}
\psfrag{xi}{$x_i$}
\psfrag{Ri}{$R_i$}
\includegraphics[width=175pt]{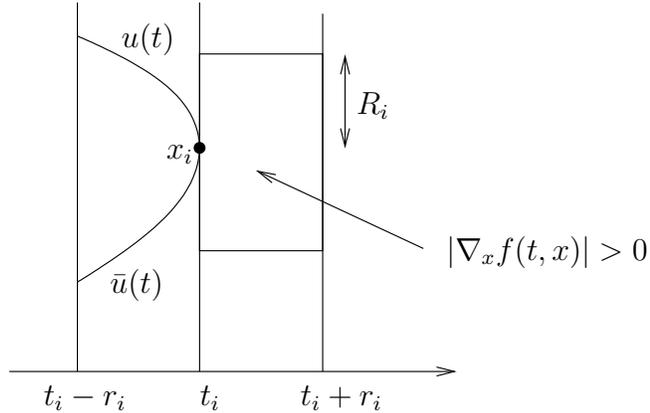}
\caption{\em The local structure of $\Gamma_f$ near $t_i$.}
\end{center}
\end{figure}
\bigskip

In the second part of this section we study some properties of the following $\e$-gradient system
\begin{equation}\llabel{ue}
\e \due(t)=-\n f(t,\ue(t)).
\end{equation}
We start by proving the existence of global solutions to Cauchy problems associated to (\ref{ue}). By {\em global} we mean a solution
defined on the whole interval $[0,T]$. 
\begin{lemma}\llabel{lm:u-global}
Under Assumption~\ref{Hp1}, for any $x\in \RR^n$ there exists a unique solution $t\mapsto \ue(t)$ to equation~(\ref{ue}), defined on the whole interval $[0,T]$, with the initial condition $\ue (0) = x$. 
Moreover, $\ue(t)$ is bounded uniformly with respect to $t$ and $\e$.
\end{lemma}
\begin{proof}
Since, by assumption, function $f$ is regular, it follows from standard arguments on ordinary differential equations that for every $\e$, the Cauchy problem associated to (\ref{ue}) has locally a unique solution $t\mapsto \ue(t)$.
Moreover, multiplying (\ref{ue}) by $\ue$ and using Assumption~\ref{Hp1} we get 
\begin{equation*}
\frac{d}{dt}|\ue(t)|^2\leq 2\frac{a_0}{\e}-2\frac{c_0}{\e}|\ue(t)|^2\,,
\end{equation*}
which in particular implies that for every $\e$ the solution $\ue$ is defined on $[0,T]$.

By a standard comparison argument it follows that
\begin{equation*}
|\ue(t)|^2\leq \frac{a_0}{c_0}+ e^{-2\frac{c_0}{\e}\,t}\big(|x|^2-\frac{a_0}{c_0}\big)\leq \max\Big\{\frac{a_0}{c_0},|x|^2\Big\}\,,
\end{equation*}
which gives the uniform boundedness of $\ue$ with respect to $t$ and $\e$. This concludes the proof.
\end{proof}

In the next proposition we deduce another important fact for the sequence $(\ue)_\e$.
\begin{proposition}\llabel{p:eu'}
Suppose that Assumption~\ref{Hp1} is satisfied. 
For every $\e$, let $\ue$ be the solution to a Cauchy problem associated to (\ref{ue}). 
Then 
\begin{equation*}
\e \dot{u}_{\e} \to 0 \quad \mbox{strongly in }L^2([0,T]),
\end{equation*}
as $\e$ goes to zero.
\end{proposition}
\begin{proof}
Let us notice that 
\begin{equation}\llabel{nfu}
-\n f(t,\ue)\due = -\frac{d}{dt}\, f(t,\ue) + f_t(t,\ue).
\end{equation}
Multiplying equation (\ref{ue}) by $\due$, integrating between $0$ and $T$, and taking into account (\ref{nfu}), we get
\begin{equation*}
\e \int_0^T |\due(t)|^2\, dt = f(0,\ue(0)) - f(T, \ue(T)) + \int_0^T f_t(t,\ue(t)) \, dt.
\end{equation*}
The conclusion follows now from the fact that the right-hand side is bounded uniformly with respect to $\e$.
\end{proof}

Now we are in a position to state the main result of this paper.
\begin{theorem}\llabel{tm:main}
Under Assumptions~\ref{Hp1}--\ref{Hp4}, let $y_0\in\RR^n$ be such that $\n f(0,y_0)=0$ and $\h f(0,y_0)$ is positive definite. 
Let $u\colon[0,T]\to\RR^n$ be the $C^2$-piecewise function given by Definition~\ref{d:u} with $u(0)=y_0$, and let $\ue\colon[0,T]\to\RR^n$ be the solution of (\ref{ue}) starting from $\ue(0)=:y_\e\in\RR^n$.
If $y_\e\to y_0$, then
\begin{equation}\llabel{tesi1}
\ue \to u \quad\mbox{ uniformly on compact sets of  }[0,T]\setminus \{t_1,\dots,t_k\}.
\end{equation}
Moreover, for every $i=1,\dots,k-1$ let $v_i$ be the heteroclinic solution of (\ref{vi})--(\ref{vii}). 
Then there exists $t_\e^i$ such that $t_\e^i\to t_i$ as $\e\to 0$ and
\begin{equation}\llabel{tesi2}
v^i_\e(s):=\ue(t_\e^i+\e s) \to v_i(s)\quad \mbox{ uniformly on compact sets of }\RR.
\end{equation}
Finally, if $\gamma_i:=v_i(\RR)\cup\{x_i,y_i\}$ represents the trajectory of $v_i$ and 
\begin{equation}\llabel{G}
G:={\rm graph }(u)\cup \bigcup_{i=1}^k (\{t_i\}\times\gamma_i),
\end{equation}
then
\begin{equation}\llabel{tesi3}
{\rm dist}((t,\ue(t)),G)\to 0 \quad\mbox{  as $\e\to 0$}\,,
\end{equation}
uniformly for $t\in[0,T]$.
\end{theorem} 
\begin{remark}\llabel{r:ilthm}
In previous theorem the following three facts are established.
First, that out of some small neighborhoods of the critical times $t_i$, the distance between the ``perturbed'' solution $\ue(t)$ and the limit function $u(t)$ is small uniformly with respect to $t$. 
Next, that in a small neighborhood of $t_i$, the solution $\ue$ belongs to a tubular neighborhood of 
the trajectory of the heteroclinic solution $v_i$.

Notice that these two facts together imply that the graph of $\ue$ approaches the completion of the graph of $u$ obtained by using the heteroclinic trajectories, defined in (\ref{G}).

The third fact is that near the critical times $t_i$, a suitable rescaled version of $\ue$ converges to the heteroclinic solution~$v_i$.
\end{remark}
\end{section}

\begin{section}{Proof of the main result }\llabel{s:result}

The proof of Theorem~\ref{tm:main} follows from some intermediate lemmas which we are going to prove.
For simplicity, we focus on the first subinterval $[0,t_1]$ and we start by showing in the next lemma
that (\ref{tesi1}) holds true.
\begin{lemma}\llabel{tm:near-u}
Under the assumptions of Theorem~\ref{tm:main}, if $\ue(0)\to u(0)$, then $\ue\to u$ uniformly on compact subsets of $[0,t_1)$.
\end{lemma}
\begin{proof}
For every $0\leq \tau<t_1$, by construction of the function $u$, there exists $\alpha=\alpha(\tau)$ such that
\begin{equation*}
\h f(t,u(t))y\cdot y\geq 2\alpha|y|^2
\end{equation*}
for every $y\in\RR^n$ and every $0\leq t\leq \tau$. (Indeed it is sufficient to take $\alpha(\tau)$ be equal to the smallest (positive) eigenvalue of $\h f(t,u(t))$ for $t\in[0,\tau]$).

By uniform continuity, there exists $\D_0>0$ such that 
\begin{equation}\llabel{ha}
\h f(t,x)y\cdot y\geq \alpha|y|^2
\end{equation}
for every $y\in\RR^n$ and every $t\in[0,\tau]$, provided that $|u(t)-x|<\D_0$.

Since $\ue(0)$ converges to $u(0)$ as $\e\to 0$, then there exists $\e_0>0$ such that $|\ue(0)-u(0)|<\D_0$, for every $\e<\e_0$. 
Let $t^*_1$ be the largest time such that $|\ue(t)-u(t)|<\D_0$ for every $t\in[0,t^*_1)$, i.e.,
\begin{equation*}
t^*_1:=\sup\{t\in[0,\tau):|\ue(t)-u(t)|<\D_0\}\,.
\end{equation*}
For every $t\in[0,t^*_1)$ and every $\e<\e_0$, subtracting $\e \dot{u}(t)$ to (\ref{ue}), we deduce that
\begin{equation*}
\e(\dot{u}_\e(t)-\dot{u}(t)) = -\n f(t,\ue(t)) + \n f(t,u(t)) -\e\dot{u}(t).
\end{equation*}
Let us multiply previous equation by $w_\e(t):=\ue(t)-u(t)$. 
By the Mean Value Theorem, and using (\ref{ha}), and the Cauchy inequality, we obtain
\begin{equation*}
\begin{array}{l}
{\displaystyle \frac{\e}{2}\frac{d}{dt}|w_\e(t)|^2 
\leq -\alpha |w_\e(t)|^2 + \frac{\e}{2}\beta + \frac{\e}{2}|w_\e(t)|^2}\,,
\end{array}
\end{equation*}
where $\beta$ is an upper bound for $|\dot{u}(t)|^2$.
Using differential inequalities, we deduce that
\begin{equation}\llabel{we}
|w_\e(t)|^2\leq \Big(|w_\e(0)|^2-\e\frac{\beta}{2\alpha-\e}\Big)e^{-(2\frac{\alpha}{\e}-1)t} + \e\frac{\beta}{2\alpha-\e}\,.
\end{equation}
It follows from (\ref{we}) that for $\e$ small enough $|\ue(t^*_1)-u(t^*_1)|<\D_0$, which, by the definition of $t^*_1$, implies that $t^*_1=\tau$.
Moreover, since by assumption $w_\e(0)\to 0$, we deduce from (\ref{we}) that $w_\e(t)\to 0$ uniformly on $[0,\tau]$ as $\e\to 0$, and this concludes the proof.
\end{proof}

In order to prove condition (\ref{tesi2}) in Theorem~\ref{tm:main}, we zoom in on a neighborhood of $t_1$ and discuss what happens.
Let $x_1$ be defined by condition~(3) in Proposition~\ref{p:u}, and let $\Lambda:=\min\{|x_1-y|: \n f(t_1,y)=0\,, y\neq x_1\}$ be the minimal distance between $x_1$ and the other critical points of $f(t_1,\cdot)$. 
Let $0<\D_1<\min\{\Lambda,R_1\}$, where $R_1$ is the constant such that inequality (\ref{ehi}) is satisfied for every $t\in {]t_1,t_1+r_1]}$, and $|x-x_1|<R_1$.

By continuity, since $u(t)$ tends to $x_1$ as $t\to t_1^-$, there exists $\bar{t}<t_1$ such that 
\begin{equation}\llabel{ud1}
|u(t)-x_1|<\frac{\D_1}{4} \quad\forall \, t\in(\bar{t},t_1)\,.
\end{equation}
Consider now an increasing sequence $\tau_h$ approaching $t_1$ to the left, with $\tau_1>\bar{t}$. 
Since, for every $h$, Lemma~\ref{tm:near-u} implies that $|\ue(t)-u(t)|\to 0$ uniformly on $[0,\tau_h]$, we deduce that there exists $\e_h>0$ such that
\begin{equation}\llabel{tauh}
|\ue(t)-u(t)|<\frac{\D_1}{4} \quad \forall \, t\in[0,\tau_h],\, \forall \,0<\e<\e_h.
\end{equation}
Let us define
\begin{equation}\llabel{tde}
t^{\D_1}_\e:= \inf\{t\geq\tau_1:|\ue(t)-x_1|\geq \D_1\},
\end{equation}
i.e., $t_\e^{\D_1}$ is the first time larger than $\tau_1$ such that $|\ue(t)-x_1|=\D_1$.
\begin{lemma}\llabel{lm:tesi2}
Let $t^{\D_1}_\e$ be defined by (\ref{tde}). Then
\begin{equation}\llabel{tde2}
t^{\D_1}_\e\to t_1 \qquad\mbox{ as }\e\to 0.
\end{equation}
\end{lemma}
\begin{proof}
We begin by proving that
\begin{equation}\llabel{tesi2-1}
\liminf_{\e\to 0}t^{\D_1}_\e\geq t_1.
\end{equation}
Let $\bar{t}<t_1$ be such that (\ref{ud1}) is satisfied.
Since by definition $\tau_1>\bar{t}$, we have that $\tau_h$ belongs to $(\bar{t},t_1)$ for every $h$. 
Then for fixed $\tau_h$ it follows from (\ref{ud1}), (\ref{tauh}), and triangular inequality, that 
\begin{equation*}
|\ue(t)-x_1|<\frac{\D_1}{2} \quad\forall\,t\in(\bar{t},\tau_h), \quad\mbox{and every $\e<\e_h$}. 
\end{equation*}
Hence, $t^{\D_1}_\e>\tau_h$, for every $0<\e<\e_h$.
Thus, $\liminf_{\e\to 0}t^{\D_1}_\e\geq \tau_h$ for every $h$, which implies~(\ref{tesi2-1}). 

On the other hand, by Proposition~\ref{p:eu'}
for a.e.\ $t^*\in[0,T]$, $|\n f(t^*,\ue(t^*))|$ tends to $0$ as $\e\to 0$ along a suitable sequence. 
In particular, this is true for a.e.\ $t^*$ in a right-neighborhood of $t_1$. 
Condition (\ref{ehi}) implies now that $|\ue(t^*)-x_1|>R_1$ for $\e$ sufficiently small. 
Let us take $\eta>0$ and choose $t^*\in{]t_1,t_1+\eta[}$.
Since $R_1>\D_1$, from the definition of $t^{\D_1}_\e$ and the regularity of $\ue$, we deduce immediately that $t^{\D_1}_\e<t^*$ for $\e$ sufficiently small.
This concludes the proof, since the result does not depend on the subsequence of $\e$ chosen.
\end{proof}

Let us observe now that for $s\in [-t^{\D_1}_\e/\e,(T-t^{\D_1}_\e)/\e]$, function 
$v^\e_1(s):=\ue(t^{\D_1}_\e+\e s)$ solves the following problem
\begin{equation}\llabel{ve1}
\begin{cases}
\dot{v}^\e_1(s)=-\n f(t^{\D_1}_\e+\e s,v^\e_1(s))\\
v^\e_1(0)=\ue(t^{\D_1}_\e).
\end{cases}
\end{equation}
Moreover, since $\ue(t^{\D_1}_\e)$ belongs to the compact set $\partial B(x_1,\D_1)$, there exists $\kappa_1\in \partial B(x_1,\D_1)$ such that, passing to a subsequence, $\ue(t^{\D_1}_\e)\to \kappa_1$ as $\e\to 0$.
Therefore, Lemma~\ref{lm:tesi2} and the Continuous Dependence Theorem imply that $v^\e_1$ converges uniformly 
on compact sets of $\RR$ to the solution $w(s)$ of the following problem:
\begin{equation}\llabel{whv}
\begin{cases}
\dot{w}(s)=-\n f(t_1,w(s))\\
w(0)=\kappa_1.
\end{cases}
\end{equation}
The next step consists in proving that $w$ is precisely (up to time-translations) the heteroclinic solution $v_1$, defined in Proposition~\ref{p:u}.
To this aim we introduce a sequence $\D_k\searrow 0$, where $\D_1$ is the constant already introduced (after the proof of Lemma~\ref{tm:near-u}), and define, for $k>1$,
\begin{equation}\llabel{tdek}
t^{\D_k}_\e := \sup\{t\leq t^{\D_1}_\e: |\ue(t)-x_1|\leq \D_k\},
\end{equation}
i.e., $t^{\D_k}_\e$ is the last time before $t^{\D_1}_\e$ such that $|\ue(t)-x_1| = \D_k$.
\begin{lemma}\llabel{lm:2}
For $s\in [-t^{\D_k}_\e/\e,(T-t^{\D_k}_\e)/\e]$, let $v^\e_k(s):= \ue(t^{\D_k}_\e+\e s)$, and let $t^{\D_1}_\e=t^{\D_k}_\e + \e S^{1,k}_\e$, for some $S^{1,k}_\e>0$.
Then $S^{1,k}_\e\to s_k<+\infty$ as $\e\to 0$ along a suitable sequence, and, for every $k$,
\begin{equation}\llabel{e:2}
v^\e_k(s)\to w(s-s_k) \qquad\mbox{ uniformly on compact subsets of }\RR\,.
\end{equation}
Moreover $s_k\to+\infty$ as $k\to\infty$. 
Finally, $\lim_{s\to-\infty}w(s)=x_1$.
\end{lemma}
\begin{proof}
We begin by observing that $t^{\D_k}_\e$ converges to $t_1$ as $\e\to 0$. 
Indeed, arguing in the same manner as for $t^{\D_1}_\e$ in the proof of Lemma~\ref{lm:tesi2}, we get that $\liminf_{\e\to 0}t^{\D_k}_\e\geq t_1$,
while, since $t^{\D_k}_\e\leq t^{\D_1}_\e$, and $t^{\D_1}_\e\to t_1$, we deduce that $\limsup_{\e\to 0}t^{\D_k}_\e\leq t_1$.

Near $(t_1,x_1)$, using the local structure of the set $\Gamma_f$ (defined in (\ref{Gf})), given by ($C_l$) and ($C_r$), we can prove that for every $k$ there exists $\eta_k>0$ such that 
\begin{equation}\llabel{lz}
\Gamma_f\cap \big([t_1-\eta_k,t_1+\eta_k]\times \overline{B}(x_1,R_1)\big)\subset [t_1-\eta_k,t_1]\times\overline{B}\big(x_1,\frac{\D_k}{2}\big)\,,
\end{equation}
where $R_1$ is the constant such that (\ref{ehi}) is satisfied for every $t\in{]t_1,t_1+r_1]}$ and $|x-x_1|<R_1$.
Next, we notice that for fixed $k$, for $\e$ sufficiently small and for $t\in [t^{\D_k}_\e, t^{\D_1}_\e]$, we have
\begin{equation*}
(t,\ue(t))\in \mathcal{S}_k:=[t_1-\eta_k,t_1+\eta_k]\times\{x\in\RR^n:\D_k\leq |x-x_1|\leq \D_1\}\,.
\end{equation*}
Moreover, we observe that $\mathcal{S}_k$ is closed and since, by (\ref{lz}), $(t,x)\in\mathcal{S}_k$ is quite distant from both $(t,u(t)$ and $(t,\overline{u}(t))$ (the two regular branches of $\Gamma_f$ near $(t_1,x_1)$), we have $\mathcal{S}_k\cap\Gamma_f=\emptyset$, so that there exists a positive constant $c_k$ such that $|\n f(t,x)|\geq c_k>0$ for every $(t,x)\in \mathcal{S}_k$.
By the fact that $(t,\ue(t))\in \mathcal{S}_k$ for $t\in [t^{\D_k}_\e, t^{\D_1}_\e]$, it follows 
\begin{equation}\llabel{nz}
|\n f(t,\ue(t))|\geq c_k>0\qquad \mbox{ for every }t^{\D_k}_\e\leq t\leq t^{\D_1}_\e\,,
\end{equation}
and for $\e$ small. 
(See also Figure~\ref{f:Sk}.)
\begin{figure}[ht]\llabel{f:Sk}
\begin{center}
\psfrag{x1}{\scriptsize$x_1$}
\psfrag{d1}{\scriptsize$\delta_1$}
\psfrag{dk}{\scriptsize$\delta_k$}
\psfrag{Gf}{\scriptsize${\displaystyle [t_1-\eta_k,t_1]\times\overline{B}(x_1,\frac{\delta_k}{2})}$}
\psfrag{Sk}{\scriptsize$\mathcal{S}_k$}
\psfrag{t1}{\scriptsize$t_1$}
\psfrag{t1ek}{\scriptsize$t_1+\eta_k$}
\psfrag{-t1ek}{\scriptsize$t_1-\eta_k$}
\includegraphics[width=300pt]{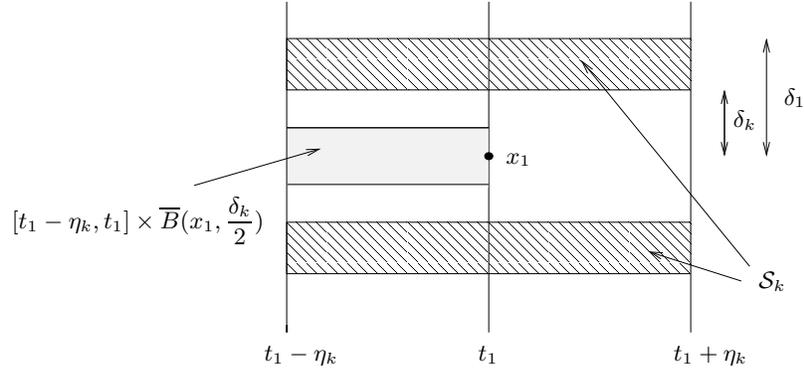}
\caption{\em The sets $\mathcal{S}_k$ and $[t_1-\eta_k,t_1+\eta_k]\times\overline{B}(x_1,R_1)$ are disjoint.}
\end{center}
\end{figure}

After these preliminaries, we prove now that $S^{1,k}_\e\to s_k<+\infty$, as $\e\to 0$ along a suitable sequence.
Indeed, for $s\in[0,1]$ let us define $w_\e(s):=\ue(s(t^{\D_1}_\e-t^{\D_k}_\e)+t^{\D_k}_\e)$.
Hence, $w_\e$ is the solution of the following problem:
\begin{equation}\llabel{Skwe}
\begin{cases}
{\displaystyle \frac{1}{S^{1,k}_\e}\, \dot{w}_\e(s)=-\n f(s(t^{\D_1}_\e-t^{\D_k}_\e)+t^{\D_k}_\e, w_\e(s))}\\
w_\e(0)=\ue(t^{\D_k}_\e).
\end{cases}
\end{equation}
Multiplying by $\dot{w}_\e$ and integrating between $0$ and $1$ we get:
\begin{equation*}
\begin{array}{c}
{\displaystyle \frac{1}{S^{1,k}_\e} \int_0^1|\dot{w}_\e(s)|^2\, ds = f(t^{\D_k}_\e, \ue(t^{\D_k}_\e))- f(t^{\D_1}_\e, \ue(t^{\D_1}_\e))+}\\
{\displaystyle \phantom{\frac{1}{S^{1,k}_\e} } +(t^{\D_1}_\e-t^{\D_k}_\e)\int_0^1 f_t(s(t^{\D_1}_\e-t^{\D_k}_\e)+t^{\D_k}_\e,w_\e(s))\, ds,}
\end{array}
\end{equation*}
where the right-hand side is bounded uniformly with respect to $\e$.
If, by contradiction, $S^{1,k}_\e\to \infty$ as $\e\to 0$ along a suitable sequence, then $(S^{1,k}_\e)^{-1}\dot{w}_\e\to 0$ strongly in $L^2(0,1)$,
which in particular implies that
\begin{equation}\llabel{Skdwe}
\frac{1}{S^{1,k}_\e}\dot{w}_\e(s)\to 0\quad\mbox{ for a.e.\ }s\in[0,1].
\end{equation}
By (\ref{nz}) and the definition of $w_\e$, for $\e$ sufficiently small, we obtain that
$$
|\n f(s(t^{\D_1}_\e-t^{\D_k}_\e)+t^{\D_k}_\e, w_\e(s))|\geq c_k>0\,,
$$ 
for every $s\in[0,1]$, which contradicts (\ref{Skwe}) and (\ref{Skdwe}).

Now we continue as for $\D_1$, and define $v^\e_{k}(s):=\ue(t^{\D_k}_\e+\e s)$, for $s\in[-t^{\D_k}_\e/\e,T-t^{\D_k}_\e/\e]$. 
It turns out that $v^\e_k$ solves
\begin{equation*}
\begin{cases}

\dot{v}^\e_{k}(s)=-\n f(t^{\D_k}_\e+\e s,v^\e_{k}(s)),\\
v^\e_{k}(0)=\ue(t^{\D_k}_\e).
\end{cases}
\end{equation*}
Since $\ue(t^{\D_k}_\e)$ belongs to the compact set $\partial B(x_1,\D_k)$, we deduce that there exists $\kappa_k\in \partial B(x_1,\D_k)$
such that $\ue(t^{\D_k}_\e)\to \kappa_k$ as $\e\to 0$ along a suitable sequence. 
It follows that, if we define $w_k$ as the solution of
\begin{equation*}
\begin{cases}
\dot{w}_k(s)=-\n f (t_1,w_k(s))\\
w_k(0)=\kappa_k,
\end{cases}
\end{equation*}
then $v^\e_{k}(s) \to w_k(s)$ uniformly on compact subsets of $\RR$, by the Continuous Dependence Theorem.
Moreover, the equality $v^\e_{k}(S^{1,k}_\e)=\ue(t^{\D_1}_\e)$ implies $w_k(s_k)=\kappa_1$. By the uniqueness of the solution of the Cauchy problem (\ref{whv}) we obtain $w_k(s)=w(s-s_k)$ and $w(-s_k)=\kappa_k$.
It follows that $w(-s_k)\to x_1$ as $k\to +\infty$.

Let now $s_\infty$ be such that $s_k\to s_\infty$, and assume by contradiction that $s_\infty<+\infty$.
Then, by continuity, $w(s_\infty)=x_1$ and, since $x_1$ is an equilibrium point, from the uniqueness it should follow $w(s)\equiv x_1$, a contradiction.

It remains to prove that $\lim_{s\to-\infty} w(s)=x_1$.
Indeed this follows from some standard facts on the $\alpha$-limit set (see, e.g., \cite{LaSalle}, \cite{Amann}).

More precisely, let $g(x):=f(t_1,x)$ for every $x\in\RR^n$, and let $E$ be the set of critical points of $g$, i.e., $E=\{x:\nabla g(x)=0\}$. 
Let us call the $\alpha$-limit set of $w$ by $\alpha(w)$.
Then $x_1\in\alpha(w)$, and we can prove that for every $y\in\alpha(w)$ we have $g(y)=g(x_1)$.
Indeed, for $y\in\alpha(w)$ there exists a sequence $\hat{s}_k$ such that $w(-\hat{s}_k)\to y$ and $\hat{s}_k\to +\infty$.
Moreover, the sequence $g(w(-\hat{s}_k))$ converges to $g(y)$ as $k\to\infty$. 
Since the map $s\mapsto g(w(s))$ is nonincreasing, there exists $a\in\RR\cup\{+\infty\} $ such that $g(w(s))\to a$, as $s\to-\infty$.
But
\begin{equation*}
\lim_{s\to-\infty} g(w(s)) = \lim_{k\to\infty} g(w(-s_k)) = \lim_{k\to\infty}g(\kappa_k)=g(x_1)\,,
\end{equation*}
so that $a=g(x_1)\in \RR$ and $a=g(y)$, for every $y\in\alpha(w)$.
It follows that $g(w(s))\leq a$ for every $s\leq 0$.
By (\ref{finf}) we deduce that the negative semiorbit is precompact. 
Then $\alpha(w)$ is connected, compact, and contained in~$E$.

Since by assumption the points of $E$ are isolated, we have $\alpha(w)=\{x_1\}$, therefore $\lim_{s\to-\infty}w(s)=x_1$.
\end{proof}
By Lemma~\ref{lm:2} there exists a subsequence $\e_k\to 0$ such that
\begin{equation*}
u_{\e_k}(t_{\e_k}^{\D_1}+\e_k s)\to w(s)\,.
\end{equation*}
Let us choose now $v_1$ satisfying (\ref{vi}) and (\ref{vii}) of Proposition~\ref{p:u}.
Then there exist $\alpha,\beta\in \RR$ such that
\begin{equation}\llabel{albe}
\{s\in\RR:v_1(s)\in\partial B(x_1,\D_1)\}\subset [\alpha,\beta]\,.
\end{equation}
{}From the fact that $w(0)=\kappa_1\in\partial B(x_1,\D_1)$ and by uniqueness, there exists unique $c\in \RR$ (defined by $v_1(c)=w(0)=\kappa_1$) such that
\begin{equation}\llabel{wv1}
w(s)=v_1(s+c)\,.
\end{equation}
In order to prove that our main result does not depend on suitable subsequences of $\e$, we will use the following lemma (given in an abstract setting).
\begin{lemma}\llabel{lm:ab}
Let $\mathcal{C}$ be a nonempty set, let $(\mathcal{F},d)$ be a metric space, and let us consider a function $G\colon {]0,\e_0[}\times\mathcal{C}\to \mathcal{F}$.
Assume that there exists $g_0\in \mathcal{F}$ such that
\begin{equation}\llabel{ab1}
\forall \e_k\to 0\; \exists \e_{k_j}\,,\exists c\in\mathcal{C}: G(\e_{k_j},c)\to g_0\,.
\end{equation}
Then for every $\e\in{]0,\e_0[}$ there exists $c_\e\in\mathcal{C}$ such that 
\begin{equation*}
G(\e,c_\e)\to g_0\,.
\end{equation*}
\end{lemma}
\begin{proof}
It is sufficient to prove that
\begin{equation}\llabel{ab0}
\inf_{c\in\mathcal{C}} d(G(\e,c),g_0)\to 0\,.
\end{equation}
Assume by contradiction that (\ref{ab0}) is not true. 
Then, there exists $\eta>0$ such that for every $k\in\NN$ there exists $\e_k\in{]0,\frac{1}{k}[}$ with
\begin{equation}\llabel{ab2}
\inf_{c\in\mathcal{C}}d(G(\e_k,c),g_0)>\eta\,.
\end{equation}
But this contradicts the assumption (\ref{ab1}), since for every subsequence $\e_{k_j}$ of $\e_k$ and for every $c\in\mathcal{C}$ we should get
$d(G(\e_{k_j},c),g_0)>\eta$.
\end{proof}
Now we are in a position to prove the main result of this paper.
\begin{proof}[Proof of Theorem~\ref{tm:main}]
Let us first concentrate on the time interval $[0,t_1]$.
Then Lemma~\ref{tm:near-u} implies that condition (\ref{tesi1}) restricted to $[0,t_1]$ is satisfied.

Let us prove now condition~(\ref{tesi3}).
For fixed $\eta\in{]0,\frac{t_1}{2}[}$ the goal is to prove that there exists $\e_0>0$ such that
\begin{equation}\llabel{tesi3p}
{\rm dist}((t,\ue(t)),G)<\eta \quad\mbox{ for every }\e<\e_0\,,
\end{equation}
uniformly with respect to $t\in[0,t_1]$. 
Indeed, let $\tau:= t_1-\frac{\eta}{2}$.
Let us take $\D_1<\eta$ and define $t_\e^{\D_1}$ as in (\ref{tde}).
We consider now $\e$ belonging to a suitable sequence tending to zero such that $\kappa_1$ is defined (as the limit of $\ue(t^{\D_1}_\e)$).
Then by Lemma~\ref{tm:near-u} there exists $\e_1>0$ such that (\ref{tesi3p}) is satisfied for every $t\in[0,\tau]$ and every $\e<\e_1$. 
Moreover, from the definition of $t_\e^{\D_1}$ we deduce that $|\ue(t)-x_1|\leq\D_1<\eta$, for every $t\in{]\tau,t_\e^{\D_1}]}$.
Hence (\ref{tesi3p}) is satisfied for every $t\in[0,t_\e^{\D_1}]$ and every $\e<\e_1$.

On the other hand, by Lemma~\ref{lm:2} and by (\ref{wv1}), there exists $\e_2>0$ such that 
\begin{equation*}
|\ue(t_\e^{\D_1}+\e s)-v_1(s+c)|<\frac{\eta}{2} \quad\mbox{ for every }0\leq s\leq S_\eta \mbox{ and }\e<\e_2\,,
\end{equation*}
where $S_\eta\in\RR$ is such that 
\begin{equation*}
|v_1(s)-y_1|<\frac{\eta}{2}\qquad\mbox{ for every }s\geq S_\eta\,.
\end{equation*}
Let $\tau_\e^1:= t_\e^{\D_1}+\e S_\eta -\e c$ and observe that it is not restrictive to assume that $\tau_\e^1>t_\e^{\D_1}$.
Hence, for $s=S_\eta$ and $\e$ small enough, 
\begin{equation*}
|\ue(\tau_\e^1)-v_1(S_\eta)|<\frac{\eta}{2}\,.
\end{equation*}
We have thus obtained that 
\begin{equation*}
{\rm dist}((t,\ue(t)),\{t_1\}\times\gamma_1)\leq \eta \qquad\mbox{ on } [t_\e^{\D_1},\tau_\e^1]\,,
\end{equation*}
recalling that $\gamma_1$ is the trajectory of $v_1$.
This, together with the fact that $\tau_\e^1\to t_1$, completes the proof of  (\ref{tesi3p}) and begins the proof of the uniform convergence in the interval $[t_1,t_2]$.
After a finite number of steps we obtain (\ref{tesi1}) and (\ref{tesi3}).

Since (\ref{tesi1}) and (\ref{tesi3}) do not depend on the particular subsequence chosen, we deduce that the result holds true for the whole sequence~$\e$.

More delicate to prove is condition (\ref{tesi2}), since the constant $c$ introduced in (\ref{wv1}) depends on the subsequence $\e_k$.
Indeed we recall that $u_{\e_k}(t_{\e_k}^{\D_1}+\e_k s)\to w(s)=v_1(s+c)$, i.e.,
\begin{equation*}
u_{\e_k}(t_{\e_k}^{\D_1}-\e_k c +\e_k s)\to v_1(s)\,.
\end{equation*}
Therefore Lemma~\ref{lm:ab} applies with $\mathcal{C}:=[\alpha,\beta]$ (see (\ref{albe})) and $\mathcal{F}$ be equal to the set of continuous functions endowed with the distance induced by the uniform convergence on compact sets, and we obtain~(\ref{tesi2}).
\end{proof}
\end{section}

\noindent {\bf Acknowledgments.}{ The author wishes to thank Gianni Dal Maso for many helpful and interesting discussions on the topic, and Massimiliano Berti for useful suggestions about the proof of Lemma~\ref{lm:evi}.}

\end{document}